\documentclass{amsart}
\begin{document}
\title[newforms]{Differential modular forms attached to newforms mod $p$}
\author{Alexandru Buium}
\def \Rp{R_p}
\def \Rpi{R_{\pi}}
\def \dpi{\d_{\pi}}
\def \bT{{\bf T}}
\def \cI{{\mathcal I}}
\def \cJ{{\mathcal J}}
\def \ZN{\bZ[1/N,\zeta_N]}
\def \tA{\tilde{A}}
\def \o{\omega}
\def \tB{\tilde{B}}
\def \tC{\tilde{C}}
\def \alph{A}
\def \bet{B}
\def \bsigma{\bar{\sigma}}
\def \y{^{\infty}}
\def \Ra{\Rightarrow}
\def \uBS{\overline{BS}}
\def \lBS{\underline{BS}}
\def \lB{\underline{B}}
\def \<{\langle}
\def \>{\rangle}
\def \hL{\hat{L}}
\def \cU{\mathcal U}
\def \cF{\mathcal F}
\def \S{\Sigma}
\def \st{\stackrel}
\def \sd{Spec_{\d}\ }
\def \pd{Proj_{\d}\ }
\def \s{\sigma_2}
\def \i{\sigma_1}
\def \bs{\bigskip}
\def \cD{\mathcal D}
\def \cC{\mathcal C}
\def \cT{\mathcal T}
\def \cK{\mathcal K}
\def \cX{\mathcal X}
\def \sX{X_{set}}
\def \cY{\mathcal Y}
\def \cS{X}
\def \cR{\mathcal R}
\def \cE{\mathcal E}
\def \tcE{\tilde{\mathcal E}}
\def \cP{\mathcal P}
\def \cA{\mathcal A}
\def \cV{\mathcal V}
\def \cM{\mathcal M}
\def \cL{\mathcal L}
\def \cN{\mathcal N}
\def \tcM{\tilde{\mathcal M}}
\def \caS{\mathcal S}
\def \cG{\mathcal G}
\def \cB{\mathcal B}
\def \tG{\tilde{G}}
\def \cF{\mathcal F}
\def \h{\hat{\ }}
\def \hp{\hat{\ }}
\def \tS{\tilde{S}}
\def \tP{\tilde{P}}
\def \tA{\tilde{A}}
\def \tX{\tilde{X}}
\def \tcS{\tilde{X}}
\def \tT{\tilde{T}}
\def \tE{\tilde{E}}
\def \tV{\tilde{V}}
\def \tC{\tilde{C}}
\def \tI{\tilde{I}}
\def \tU{\tilde{U}}
\def \tG{\tilde{G}}
\def \tu{\tilde{u}}
\def \chu{\check{u}}
\def \tx{\tilde{x}}
\def \tL{\tilde{L}}
\def \tY{\tilde{Y}}
\def \d{\delta}
\def \e{\chi}
\def \bW{\mathbb W}
\def \bV{{\mathbb V}}
\def \bF{{\bf F}}
\def \bE{{\bf E}}
\def \bC{{\bf C}}
\def \bO{{\bf O}}
\def \bR{{\bf R}}
\def \bA{{\bf A}}
\def \bB{{\bf B}}
\def \cO{\mathcal O}
\def \ra{\rightarrow}
\def \bx{{\bf x}}
\def \f{{\bf f}}
\def \bX{{\bf X}}
\def \bH{{\bf H}}
\def \bS{{\bf S}}
\def \bF{{\bf F}}
\def \bN{{\bf N}}
\def \bK{{\bf K}}
\def \bE{{\bf E}}
\def \bB{{\bf B}}
\def \bQ{{\bf Q}}
\def \bd{{\bf d}}
\def \bY{{\bf Y}}
\def \bU{{\bf U}}
\def \bL{{\bf L}}
\def \bQ{{\bf Q}}
\def \bP{{\bf P}}
\def \bR{{\bf R}}
\def \bC{{\bf C}}
\def \bD{{\bf D}}
\def \bM{{\bf M}}
\def \bZ{{\mathbb Z}}
\def \xtoleqr{x^{(\leq r)}}
\def \hU{\hat{U}}
\def \k{\kappa}
\def \ee{\overline{p^{\k}}}

\newtheorem{THM}{{\!}}[section]
\newtheorem{THMX}{{\!}}
\renewcommand{\theTHMX}{}
\newtheorem{theorem}{Theorem}[section]
\newtheorem{corollary}[theorem]{Corollary}
\newtheorem{lemma}[theorem]{Lemma}
\newtheorem{proposition}[theorem]{Proposition}
\theoremstyle{definition}
\newtheorem{definition}[theorem]{Definition}
\theoremstyle{remark}
\newtheorem{remark}[theorem]{Remark}
\newtheorem{example}[theorem]{\bf Example}
\numberwithin{equation}{section}
\address{University of New Mexico \\ Albuquerque, NM 87131}
\email{buium@math.unm.edu, arnab@math.unm.edu} 
\maketitle

\begin{abstract}
In a previous paper \cite{eigen} we attached to  classical  complex newforms  $f$ of weight $2$  certain $\d_p$-modular forms $f^{\sharp}$ (in the sense of \cite{difmod, book})   of  order $2$ and weight $0$; the forms $f^{\sharp}$ can be viewed as ``dual" to $f$ and played a key role in some of the applications of the theory \cite{local}.  The aim of this paper is to provide a higher weight version of this ``$\sharp$-duality," by attaching  to  classical newforms mod $p$, $\overline{f}$, of weight $\kappa$ between $3$ and $p$,  $\d_{\pi}$-modular forms $f^{\sharp}$ of order $2$ and weight $-\kappa'$, with $\kappa'$ between $1$ and $p-2$. \end{abstract}

\section{Introduction and main result}

In a series of papers beginning with \cite{char} an arithmetic analogue of  differential equations was developed in which derivatives of functions were replaced by  Fermat quotient operators $\d_{\pi}$; the latter act on   complete discrete valuation rings $R_{\pi}$ of unequal characteristic and uniformizer $\pi$, by the formula $\d_{\pi}x:=\frac{\phi(x)-x^p}{\pi}$, where $\phi:R_{\pi}\ra R_{\pi}$ is a ring endomorphism whose reduction mod $\pi$ is the $p$-power Frobenius.
In particular spaces of $\d_{\pi}$-modular forms
 were introduced and studied (cf. \cite{difmod, book} for $\pi=p$ and \cite{over} for arbitrary $\pi$).

A basic construction in \cite{eigen} attaches to  classical newforms $f=\sum a_n q^n$ over ${\mathbb C}$, of level $N$ (i.e. on $\Gamma_1(N)$) and weight $2$,  some $\d_p$-modular forms $f^{\sharp}$ of level $N$, order $2$, and weight $0$, with $\d_p$-Fourier expansion nicely expressible in terms of $f^{(-1)}(q):=\sum \frac{a_n}{n}q^n$. The forms $f^{\sharp}$ in \cite{eigen} are arithmetic-differential objects that have no classical analogue but, rather, can be viewed as ``dual" to the  classical objects $f$; by the way, the forms $f^{\sharp}$ introduced in \cite{eigen}  played a key role in \cite{local} where  the theory of arithmetic differential equations was used to prove   finiteness results for Heegner-like points.

The aim of this paper is to extend this ``$\sharp$-duality" to the higher weight case.  
Specifically let $N$ be an integer with $(N,p)=1, N>4$, let $k$ be an algebraic closure of ${\mathbb F}_p$,  and let $\overline{f}$ be a newform of level $N$, over $k$, of  weight $\kappa$,  with $3\leq \kappa \leq p$. 
Then,  by a construction due to Serre \cite{Gross},  there exists   a newform $f$ of level $Np$ over ${\mathbb C}$, of weight $2$,  lifting $\overline{f}$, whose character is compatible, in a sense to be reviewed later, with the Teichm\"{u}ller character. (The compatibility requires an embedding of the field $K_f$ generated by the Fourier coefficients of $f$ into some $K_{\pi}:=R_{\pi}[1/p]$; the considerations below will be valid for a suitable embedding of the Galois closure of $K_f$ into a suitable $K_{\pi}$.) Let us say that $f$ is a {\it Serre lift} of $\overline{f}$.
Also let us say that an integer $\kappa'$ is a {\it conjugate} of an integer  $3\leq \kappa\leq p$ if 
$1\leq \kappa'\leq p-2$ and there exists an integer $c$ coprime to $p-1$ such that  
$\kappa'\equiv c(2-\kappa)$ mod $p-1$.
We recall from \cite{book,over} that if $q,\d_{\pi}q,...,\d_{\pi}^nq$ are variables then any $\d_{\pi}$-modular form $g$ possesses a $\d_{\pi}$-Fourier expansion 
$E(g)$ in $R_{\pi}((q))[\d_{\pi}q,...,\d_{\pi}^nq]\h$ (where $\h$ always means here $p$-adic completion). On the other hand $\phi$ acts on $K_{\pi}[[q]]^{\infty}:=\cup_n K_{\pi}[[q,\d_{\pi}q,...,\d_{\pi}^nq]]$  by $\phi(q)=q^p+\pi\d_{\pi}q$, $\phi(\d_{\pi}q)=(\d_{\pi}q)^p+\pi\d_{\pi}^2q$, etc. In particular $K_{\pi}[[q]]^{\infty}$ has a natural structure of left module over the polynomial ring $K_{\pi}[\phi]$. So in particular, for any polynomial $P(\phi)\in K_{\pi}[\phi]$ and any series $g(q)\in K_{\pi}[[q]]$ (e.g., for $g=f^{(-1)}(q)$) it makes sense to consider the series $P(\phi)g(q)\in K_{\pi}[[q]]^{\infty}$.
The above concepts will be reviewed in detail in the body of the paper. 
We will then prove the following:

\begin{theorem}\label{main}
For any newform $\overline{f}$ over $k$,  of level $N$ and weight $\kappa$,  with $3\leq \kappa \leq p$, there exists 
a non-zero $\d_{\pi}$-modular form $f^{\sharp}$ of  level $N$, order $2$, and weight
$-\kappa'$,   with $\kappa'$  a conjugate of $\kappa$,  such that $f^{\sharp}$ has  a $\d_{\pi}$-Fourier expansion of the form 
$$E(f^{\sharp})=\sum_{\sigma}P_{\sigma}(\phi)(f^{(-1)}(q))^{\sigma},$$
with $P_{\sigma}(\phi)$   polynomials in $K_{\pi}[\phi]$, $f$  a Serre lift of $\overline{f}$, and $\sigma$ running through the set of all embeddings of $K_f$ into ${\mathbb C}$.
\end{theorem}

The strategy of proof of our Theorem is as follows. The Serre lift $f$ being of weight $2$ and level $Np$ one can use the method in  \cite{eigen} to attach to it an element $f^{\sharp}$ in the ring of {\it $\d_{\pi}$-modular forms  of level $Np$ and weight $0$}. Using a construction from \cite{over} one can interpret this $f^{\sharp}$ as a sum of $\d_{\pi}$-modular forms of level $N$ and weights $0,-1,-2,...,-p+2$. Finally, an analysis of the action of the diamond operators shows that  this sum of $\d_{\pi}$-modular forms reduces to one term only; that term turns out to have weight $-\kappa'$ with $\kappa'$ conjugate to $\kappa$. Note that each given $\kappa$ has a priori several conjugates and the theorem does not tell which of these conjugates gives the weight of $f^{\sharp}$; be that as it may the weight $-\kappa'$ of $f^{\sharp}$ is always non-zero (in contrast to the weight of $f^{\sharp}$ in \cite{eigen,over} which is equal to $0$). 
Finally note that, as in \cite{eigen},  $f^{\sharp}$ in our Theorem \ref{main} is not unique. 

The plan of the paper is as follows. In section 2 we review  $\pi$-jet spaces and $\d_{\pi}$-modular forms following \cite{char,over}. In section 3 we review (and complement) some   (complex and arithmetic) facts  about classical modular forms following \cite{Gross, over}. In section 4 we construct our forms $f^{\sharp}$. In section 5  we use some of the tools developed in the paper to construct a remarkable homomorphism. Indeed, we will recall, in that section, the definition of the ring $S_{\heartsuit}^{\infty}$ of {\it Igusa $\d_{\pi}$-modular functions of level $N$} introduced in \cite{igusa} and we will prove:

\begin{theorem}\label{curious}
There is a natural  homomorphism from the ring of $\d_{\pi}$-modular forms of level $Np$ and weight $0$ to the  ring 
$S_{\heartsuit}^{\infty}\otimes_{R_p}R_{\pi}$. This homomorphism commutes with $\d_{\pi}$ and  raises orders by $1$.\end{theorem}

Cf. Theorem \ref{natur} for a more precise statement.

\bigskip

{\bf Acknowledgement}.
The  author would like to acknowledge partial support from the Simons Foundation
(award 311773),  from the Institut des Hautes Etudes Scientifiques in Bures sur Yvette,   and from the Romanian National Authority
for Scientific Research (CNCS - UEFISCDI, project number
PN-II-ID-PCE-2012-4-0201).

\section{Review of differential modular forms}
In this section  we review the basic definitions of $\pi$-jet spaces and $\d_{\pi}$-modular forms introduced in  \cite{char, difmod, over}.
Throughout  this paper $p\geq 5$ is a fixed prime.
For any ring $A$ we denote the $p$-adic completion of $A$ by $\widehat{A}$.
We denote by $R_p=\widehat{\mathbb Z_p^{ur}}$ the $p$-adic completion of the maximum unramified extension of ${\mathbb Z}_p$.
We set $K_p=R_p[1/p]$ (fraction field of $R_p$) and $k=R_p/pR_p$ (residue field of $R_p$); so $k$ is  an algebraic closure of ${\mathbb F}_p$. 
 
 We need a ramified version of the above (which is slightly more general than the one encountered in \cite{char,over}). Indeed let $F/{\mathbb Q}$ be a normal finite extension with ring of integers $\cO_F$ and let $\wp\subset \cO_F$ be a prime dividing $p$. There exists a (not necessarily unique) ring automorphism $\phi=\phi_{\wp}$ of $\cO_F$ preserving $\wp$ and inducing the $p$-power Frobenius on $\cO_F/\wp$. Let $\cO_{\wp}$ be the completion of the localization of $\cO$ at $\wp$ and $\pi\in \cO_F$ a uniformizer of $\cO_{\wp}$. 
 Then $\phi:\cO_F\ra \cO_F$ induces an automorphism of $\cO_{\wp}$ whose reduction mod $\pi$ is the $p$-power Frobenius on $\cO_{\wp}/(\pi)$. For any discrete valuation ring $\cO'$ which is a finite unramified extension of $\cO_{\wp}$, $\cO'$ is \'{e}tale over $\cO_{\wp}$ so $\phi:\cO_{\wp}\ra \cO_{\wp}$ lifts uniquely to a ring automorphism $\phi:\cO'\ra \cO'$ whose reduction mod $\pi$ is the $p$-power Frobenius. Taking limit and completion we get an automorphism $\phi$ of the ring $\widehat{\cO_{\wp}^{ur}}$, completion of the maximum unramified extension of $\cO_{\wp}$, whose reduction mod $\pi$ is the $p$-power Frobenius on the residue field. 
 Let $K_{\pi}=K_p(\pi)$ and let $R_{\pi}$ be the valuation ring of $K_{\pi}$. Since $R_{\pi}\subset \widehat{\cO_{\wp}^{ur}}$ is an unramified extension of complete discrete valuation rings with the same residue field it follows that $R_{\pi}= \widehat{\cO_{\wp}^{ur}}$ so $\phi$ is an autmorphism of $R_{\pi}$. We let $e$ be the ramification index of $R_p\subset R_{\pi}$. Note that the resulting embedding $\rho:\cO_F\ra R_{\pi}$ depends on $F$ and $\wp$ only. 
 
  Consider the polynomial
 $$C_{\pi}(X,Y):=\pi^{-1}(X^p+Y^p-(X+Y)^p)\in \bZ[p/\pi][X,Y].$$
A $\pi$-{\it derivation} from an ${\mathbb Z}[p/\pi]$-algebra
 $A$ into an $A-$algebra $B$  is a map $\d_{\pi}:A \ra B$ such that $\d_{\pi}(1)=0$ and
\[\begin{array}{rcl}
\d_{\pi}(x+y) & = &  \d_{\pi} x + \d_{\pi} y
+C_{\pi}(x,y)\\
\d_{\pi}(xy) & = & x^p \cdot \d_{\pi} y +y^p \cdot \d_{\pi} x
+\pi \cdot \d_{\pi} x \cdot \d_{\pi} y,
\end{array}\] for all $x,y \in A$. For a
$\pi-$derivation we  denote by $\phi:A \ra B$ the map
$\phi(x)=x^p+\pi \d_{\pi} x$; then $\phi$ is a ring homomorphism. A
$\d_{\pi}$-{\it prolongation sequence} is a sequence $S^*=(S^n)_{n \geq 0}$ of  ${\mathbb Z}[p/\pi]$-algebras $S^n$, $n
\geq 0$, together with ${\mathbb Z}[p/\pi]$-algebra homomorphisms  $\varphi:S^n \ra
S^{n+1}$ and $\pi-$derivations $\d_{\pi}:S^n \ra S^{n+1}$ (where $S^{n+1}$ is viewed as an $S^n$-algebra via $\varphi$) such that
$\d_{\pi} \circ \varphi=\varphi \circ \d_{\pi}$ on $S^n$ for all $n$.  A morphism
of $\d_{\pi}$-prolongation sequences, $u^*:S^* \ra \tilde{S}^*$ is a sequence
$u^n:S^n \ra \tilde{S}^n$ of ${\mathbb Z}[p/\pi]$-algebra
 homomorphisms such that $\delta_{\pi}
\circ u^n=u^{n+1} \circ \d_{\pi}$ and $\varphi \circ u^n=u^{n+1} \circ
\varphi$. Let $W$ be the ring of polynomials $\bZ[\phi]$ in the
indeterminate $\phi$.  For $w=\sum_{i=0}^r a_i \phi^i\in W$  (respectively for $w$ with $a_i \geq 0$), $S^*$ a $\d_{\pi}$-prolongation
sequence, and $x \in (S^0)^{\times}$ (respectively $x \in S^0$) we
can consider the element $x^w:=\prod_{i=0}^r \varphi^{r-i}
\phi^i(x)^{a_i} \in (S^r)^{\times}$ (respectively $x^w \in S^r$).
On the other hand we may consider the
$\pi-$derivation $\d_{\pi}:R_{\pi} \ra R_{\pi}$ given by $\d_{\pi} x=(\phi(x)-x^p)/\pi$. One can consider the  $\d_{\pi}$-prolongation sequence $R_{\pi}^*$ where
$R_{\pi}^n=R_{\pi}$ for all $n$. By a $\d_{\pi}$-{\it prolongation sequence over $R_{\pi}$} we
understand a prolongation sequence $S^*$ equipped with a morphism
$R_{\pi}^* \ra S^*$. From now on
all our $\d_{\pi}$-prolongation sequences are
assumed to be over $R_{\pi}$.
For any affine $R_{\pi}$-scheme of finite type $Y=Spec\ A$ there exists
a (unique) $\d_{\pi}$-prolongation sequence, $A^*=(A^n)_{n \geq 0}$, with $A^0=A$ such that for any $\d_{\pi}$-prolongation sequence $B^*$
 and any $R_{\pi}$-algebra homomorphism $u:A \ra B^0$ there exists a unique morphism of $\d_{\pi}$-prolongation sequences $u^*:A^* \ra B^*$ with $u^0=u$.
We define the $\pi$-{\it jet spaces} $J^n_{\pi}(Y)$ of $Y$ as the formal schemes
$J^n_{\pi}(Y):=Spf\ \hat{A^n}$. This construction globalizes to the case when $Y$ is not necessarily affine (such that the construction commutes, in the obvious sense, with open immersions). A $\d_{\pi}$-{\it function} on $Y$ (of order $n$) is 
an element $f\in \cO(J^n_{\pi}(Y))$, equivalently 
a morphism
$f:J^n_{\pi}(Y)\ra \widehat{{\mathbb A}^1}$ over $R_{\pi}$; such a map induces a map $f_*:Y(R_{\pi})\ra R_{\pi}$. (Indeed any $R_{\pi}$-point of $Y$ canonically lifts, by the universality property, to an $R_{\pi}$-point of $J^n_{\pi}(Y)$ which can then be mapped by $f_*$ to an $R_{\pi}$-point of the affine line.) If $Y$ is smooth $f_*$ uniquely determines $f$; we then write $f_*=f$ and we also refer to $f_*$ as a $\d_{\pi}$-{\it function}.
 If $Y=G$ is a smooth group scheme over $R_{\pi}$ then 
  a $\d_{\pi}$-{\it character} of $G$ is a homomorphism
 $f:J^n_{\pi}(G)\ra \widehat{{\mathbb G}_a}$; 
 such an $f$ induces a group homomorphism $f=f_*:G(R_{\pi})\ra R_{\pi}$, still referred to as a $\d_{\pi}$-character.
 Note that the $\pi$-jet spaces $J^n_{\pi}(Y)$
only depend on the $\pi$-adic completion of $Y$ and not on $Y$ so one can introduce $\pi$-jet spaces $J^n_{\pi}({\mathcal Y})$ attached
to formal $\pi$-adic schemes ${\mathcal Y}$ over $R_{\pi}$ which are locally $\pi$-adic completions of schemes of finite type over $R_{\pi}$.

 Next we review differential modular forms ($\d_p$-modular forms and $\d_{\pi}$-modular forms) following \cite{difmod,over}. Recall the modular curve $X_1(N)_{R_p}$ over $R_p$ with $(N,p)=1$, $N>4$; it is a smooth curve and it carries a line bundle $L$ such that the spaces of sections of $L^{\kappa}$ identify with the spaces of modular forms on $\Gamma_1(N)$ defined over $R_p$ of weight $\kappa$; cf. \cite{Gross}, p. 450; in loc. cit. $L$ was denoted by $\omega$.
The curve $X_1(N)_{R_p}$ contains two (disjoint) closed subsets of interest: the {\it cusp locus}  $(cusps)$ and the {\it supersingular locus} $(ss)$.
On $Y_1(N)=X_1(N)\backslash (cusps)$ the line bundle $L$ identifies with $u_*\Omega^1_{E/Y_1(N)}$  where $u:E\ra Y_1(N)$ is the corresponding universal elliptic curve.
Next consider
  an affine open set $X \subset X_1(N)_{R_p}$, disjoint from  the supersingular locus (but not necessarily from the cusps!),  and consider the restriction of $L$ to $X$ which we continue to denote by $L$. We can  consider
the affine $X$-scheme $V:=Spec\left(\bigoplus_{n \in {\mathbb Z}}L^{\otimes n}\right)\ra X$; then $V$ is naturally a principal (Zariski locally trivial)  ${\mathbb G}_m$-bundle. A $\d_{\pi}$-{\it modular function} on $X$, of level $N$ and order $r$ is a $\d_{\pi}$-function $V(R_{\pi})\ra R_{\pi}$ of order $r$. We say that  a $\d_{\pi}$-modular function $f$ has {\it weight}  $w\in W$ 
 if 
 $\lambda \star f=\lambda^w f$ for $\lambda\in R_{\pi}^{\times}$ 
 where $\star$ is the natural ${\mathbb G}_m$-action on $V$;
   $\d_{\pi}$-modular functions possessing weights are called $\d_{\pi}$-{\it modular forms}. We denote by $M^r_{\pi}(w)$ the spaces of  $\d_{\pi}$-modular forms of weight $w$ and order $r$.
Recall from \cite{over} that, for  variables denoted by
$q, \d_{\pi}q,...,\d_{\pi}^rq$, 
 one has natural injective $\d_{\pi}$-Fourier expansion maps
$$E:M^r(w)\ra R_{\pi}((q))^r:=R_{\pi}((q))[\d_{\pi}q,...,\d_{\pi}^rq]\h,$$
induced by the universality property of $(\cO(J^r_{\pi}(V)))_{r\geq 0}$ and the obvious structure of prolongation sequence of $(R_{\pi}((q))^r)_{r\geq 0}$.

\section{Review of classical modular forms}

In this section we review (and complement) some basic facts about classical modular forms in \cite{DI, Gross, over}.

\subsection{Complex modular forms \cite{DI}} We denote by $M_m(\Gamma_1(M),{\mathbb C})$ the space of classical complex modular forms of weight $m$ on $\Gamma_1(M)$ and by $T_m(n)$ and $\langle n \rangle$, $n\geq 1$, the Hecke operators and the diamond operators acting on this space. The diamond operators belong to the ring generated by the Hecke operators. An eigenform is an element $f\in M_m(\Gamma_1(M),{\mathbb C})$ which is an eigenvector for all $T_m(n)$ (and hence also for all $\langle n \rangle$).  Let $S_m(\Gamma_1(M),{\mathbb C})$ be the space of cusp forms. An eigenform has $a_1\neq 0$.  
If a cuspidal eigenform $f=\sum a_n q^n$ is normalized by $a_1=1$ then 
$T_m(n)f=a_nf$ for $n\geq 1$  and there exists a Dirichlet character $\epsilon$ mod $M$ such that $\langle n \rangle f =\epsilon(n) f$ for $n\geq 1$. Let $f$ be a normalized  cuspidal eigenform as above; one lets $\cO_f\subset {\mathbb C}$ denote the ring generated by all $a_n$ (it is finitely generated) and one lets $K_f$ be the fraction field of $\cO_f$. The Dirichlet character $\epsilon$ attached to $f$  takes then values in $\cO_f$.  A newform is a cuspidal normalized eigenform such that the system $\{a_l;l\not|M\}$ does not occur for a cuspidal normalized eigenform  form of the same weight and smaller level dividing $M$.  
Let $X_1(M)$ be the complete complex modular curve for $\Gamma_1(M)$; as an algebraic curve it is defined over ${\mathbb Q}$ hence its Jacobian $J_1(M)$ is also defined over ${\mathbb Q}$. The diamond operators act on $X_1(M)$ and are defined over ${\mathbb Q}$ and hence induce endomorphisms of $J_1(M)$ defined over ${\mathbb Q}$.
We recall the following basic:

\begin{theorem}\label{eichshim}
(Eichler-Shimura) Let $f\in S_2(\Gamma_1(M),{\mathbb C})$  be a newform with character $\epsilon$. Then there exists an abelian variety $A_{f/{\mathbb Q}}$ over ${\mathbb Q}$, a ring homomorphism $$\iota:\cO_f\ra End(A_{f/{\mathbb Q}}/{\mathbb Q}),$$  and a dominant homomorphism $$\pi:J_1(M)_{\mathbb Q}\ra A_{f/{\mathbb Q}},$$ defined over ${\mathbb Q}$, such that for any prime $l$ we have  commutative diagrams
$$\begin{array}{rcl}
J_1(M)_{\mathbb Q} & \stackrel{T(l)_*}{\longrightarrow} & J_1(M)_{\mathbb Q}\\
\pi \downarrow & \  & \downarrow \pi\\
A_{f/{\mathbb Q}} & \stackrel{\iota(a_l)}{\longrightarrow} & A_{f/{\mathbb Q}}
\end{array}
\begin{array}{rcl}
J_1(M)_{\mathbb Q} & \stackrel{\langle l \rangle_*}{\longrightarrow} & J_1(M)_{\mathbb Q}\\
\pi \downarrow & \  & \downarrow \pi\\
A_{f/{\mathbb Q}} & \stackrel{\iota(\epsilon(l))}{\longrightarrow} & A_{f/{\mathbb Q}}
\end{array}
$$
Moreover the image of the pull-back map
\[\pi^*:H^0(A_{f/\bC},\Omega) \ra
H^0(J_1(N)_{\bC},\Omega) \simeq S_2(\Gamma_1(N),\bC)\] is the
$\bC-$linear span of $\{f^{\sigma}\ |\ \sigma:K_f \ra \bC\}$.
\end{theorem}

Here $T(l)_*,\langle l \rangle_*\in End(J_1(M)_{\mathbb Q}/{\mathbb  Q})$ are the naturally induced Hecke and diamond  operators.

\subsection{Model  of $X_1(Np)$} 
For  modular curves and forms we use the standard notation from \cite{Gross}. Let $N>4$, $(N,p)=1$.
Recall  that the modular curve  $X_1(Np)$ over ${\mathbb C}$ has a model (still denoted by $X_1(Np)$ in what follows) over ${\mathbb Z}[1/N,\zeta_p]$
considered in  \cite{Gross}, p. 470.
Recall some of the main properties of $X_1(Np)$.
First $X_1(Np)$ is a regular scheme proper and flat  of relative dimension $1$ over ${\mathbb Z}[1/N,\zeta_p]$ 
(where $\zeta_p$ is a primitive $p$-th root of unity)
and is smooth over ${\mathbb Z}[1/Np,\zeta_p]$. Also the special fiber of $X_1(Np)$ at $p$  is a union of two smooth projective curves $I$ and $I'$ crossing transversally at a finite set $\Sigma$ of points. Furthermore $I$ is isomorphic to the Igusa curve $I_1(N)$ in \cite{Gross}, p. 160, so $I$ is the smooth compactification of the curve classifying triples $(E,\alpha, \beta)$ with $E$ an elliptic curve over a scheme of characteristic $p$, and $\alpha:\mu_N \ra E$, $\beta:\mu_p \ra E$ are embeddings (of group schemes). Similarly  $I'$ is the smooth compactification of the curve classifying triples $(E,\alpha, b)$ with $E$ an elliptic curve over a scheme of characteristic $p$, and $\alpha:\mu_N \ra E$, $b:{\mathbb Z}/p{\mathbb Z}\ra E$  are embeddings.
Finally $\Sigma$ corresponds to the supersingular locus on the corresponding curves.

Let now $F_0={\mathbb Q}(\zeta_p)$  and ${\pi_0}=1-\zeta_p$ and consider the   embedding of ${\mathbb Z}[\zeta_N, \zeta_p, 1/N]$ into $R_{{\pi_0}}$ (hence of ${\mathbb Z}[\zeta_N,1/N]$ into $R_p$). 
 Denote by $X_1(Np)_{R_{{\pi_0}}}$ the base change of $X_1(Np)$ to $R_{{\pi_0}}$.
Recall that the modular curve $X_1(N)$  over ${\mathbb C}$ has a natural smooth projective model (still denoted by $X_1(N)$)
over ${\mathbb Z}[1/N]$ such that
$$Y_1(N):=X_1(N)\backslash (cusps)$$
 parameterizes pairs $(E,\alpha)$
consisting of elliptic curves $E$ with an embedding $\alpha:\mu_N \ra E$. The morphism $X_1(Np)\ra X_1(N)$ over ${\mathbb C}$  induces a
morphism
$$e: X_1(Np)_{R_{{\pi_0}}}\backslash \Sigma \ra X_1(N)_{R_{{\pi_0}}}\backslash (ss)$$
 over $R_{{\pi_0}}$, where $(ss)$ is the supersingular locus in the closed fiber of $X_1(N)_{R_{{\pi_0}}}$. (See \cite{over} for details.)  Let $X \subset X_1(N)_{R_p} \backslash (ss)$ be an affine open set,
 $X_{R_{{\pi_0}}}:=X \otimes_{R_p} R_{{\pi_0}}\subset  X_1(N)_{R_{{\pi_0}}} \backslash (ss)$ its base change to $R_{{\pi_0}}$,
 and $X_!:=e^{-1}(X_{R_{{\pi_0}}})$. Denote by ${\mathcal X}_{R_{{\pi_0}}}$  the ${\pi_0}$-adic completion of $X_{R_{{\pi_0}}}$.
Also note that the ${\pi_0}$-adic completion of $X_!$ has two connected components; let ${\mathcal X}_!$ be the component whose reduction mod ${\pi_0}$ is contained in $I \backslash \Sigma$. We get a  morphism $e:{\mathcal X}_! \ra {\mathcal X}_{R_{{\pi_0}}}$.

 Now the diamond operators $\langle n \rangle$ on $X_1(Np)_{\mathbb Q}$ 
 give rise
 to diamond operators  $\langle n \rangle_N$ (for $n\in (\bZ/N\bZ)^{\times}$) and $\langle d \rangle_p$ (for $d \in (\bZ/p\bZ)^{\times}$) on $X_1(Np)_{\mathbb Q}$  corresponding to the Chinese remainder theorem spitting 
 $(\bZ/Np\bZ)^{\times} \simeq (\bZ/N\bZ)^{\times} \times (\bZ/p\bZ)^{\times}$.
 Recall from \cite{Gross} that 
 $G:=({\mathbb Z}/p{\mathbb Z})^{\times}$ acts via $\langle\ \rangle_p$ on the covering $X_1(Np) \ra X_1(N)$ over ${\mathbb Z}[1/N,\zeta_p]$; this action preserves the Igusa curve $I$ and induces on $I$ the usual diamond operators.
 So $G$ acts on the covering $e:{\mathcal X}_!\ra {\mathcal X}_{R_{{\pi_0}}}$.   
 Also, as a matter of notation, if $\epsilon:(\bZ/Np\bZ)^{\times}\ra {\mathbb C}$ is a Dirichlet character we denote by $\epsilon_N:(\bZ/N\bZ)^{\times}\ra {\mathbb C}$ and $\epsilon_p:(\bZ/p\bZ)^{\times}\ra {\mathbb C}$ the induced Dirichlet characters.

\subsection{Model of $J_1(Np)$}
In the discussion below we continue to assume that ${\pi_0}=1-\zeta_p$ and we continue to consider the embedding  embedding ${\mathbb Z}[\zeta_{N},\zeta_{p},1/N]\ra  R_{{\pi_0}}$.
 Let $A_{K_{{\pi_0}}}=J_1(Np)_{K_{{\pi_0}}}$ be the Jacobian of $X_1(Np)_{K_{{\pi_0}}}$ over $K_{{\pi_0}}$. Let $A_{R_{{\pi_0}}}$ be the N\'{e}ron model of $A_{K_{{\pi_0}}}$ over $R_{{\pi_0}}$; cf. \cite{BLR}, p. 12.
 By the N\'{e}ron property the Abel-Jacobi map $X_1(Np)_{R_{{\pi_0}}}\ra J_1(Np)_{R_{{\pi_0}}}$ that sends $\infty$ into $0$ can be extended to a morphism 
 $$X_1(Np)_{R_{{\pi_0}}}\backslash \Sigma\ra A_{R_{{\pi_0}}}$$
 and hence induces a morphism
 $${\mathcal X}_!\ra {\mathcal A}:=(A^0_{R_{{\pi_0}}})\hat{\ }$$
 where $A^0_{R_{{\pi_0}}}$ is the connected component of $A_{R_{{\pi_0}}}$.
 Note that, since ${\pi_0}=1-\zeta_p$, one has by \cite{BLR}, pp. 246 and 267, that $A^0_{R_{{\pi_0}}}$ coincides with $Pic^0_{X_1(Np)/R_{{\pi_0}}}$ and is a semi-abelian scheme
 whose special fiber modulo its maximum torus is a product of the Jacobians of $I$ and $I'$ respectively.
 
 \subsection{Pulling back to  arbitrary $R_{\pi}$.}
 Consider an arbitrary number field $F$ as in section 2, containing $F_0={\mathbb Q}(\zeta_p)$, and let $\pi$ be as in section 2. We will allow ourselves to enlarge $F$ and change $\pi$ repeatedly in what follows. 
 Denote by
 \begin{equation}
 \label{doarme}
 X_{R_{\pi}},X_!,\cX_!,A_{R_{\pi}},A^0_{R_{\pi}},\cA\end{equation}
 the objects obtained from the corresponding objects over $R_{\pi_0}$ by base change
  via $R_{\pi_0}\subset R_{\pi}$. Note that the latter extension of rings can be ramified so $A_{R_{\pi}}$ is not necessarily the N\'{e}ron model of $A_{K_{\pi}}$; but all the schemes (respectively formal schemes) above continue to be smooth over $R_{\pi}$ and $A_{R_{\pi}}^0$ is semiabelian, equal to $Pic^0_{X_1(Np)/R_{\pi}}$.
 
 \subsection{Schemes attached to a fixed complex newform}\label{model}
 Let $f\in S_2(\Gamma_1(Np),{\mathbb C})$ be a newform and assume we are given an embedding $\rho:\cO_f \ra R_{\pi}$ with $\pi$ as above. Let $A_{K_{\pi}}\ra A_{f/K_{\pi}}$ be the induced morphism of abelian varieties over $K_{\pi}$ obtained by base change from the morphism of abelian varieties over ${\mathbb Q}$, $J_1(Np)_{\mathbb Q}\ra A_{f/{\mathbb Q}}$.
   Let $A_{f/R_{\pi}}$ be the N\'{e}ron model of 
 $A_{f/K_{\pi}}$, let $A^{0}_{f/R_{\pi}}$ be the connected component of the latter
 and let $A^0_{R_{\pi}}\ra A^0_{f/R_{\pi}}$ be the induced morphism.
 By the theory of Raynaud extensions \cite{FC}, pp. 33-34, there exists a group scheme  $C$ over $R_{\pi}$, which is an extension of an abelian scheme $B$ over $R_{\pi}$ by a torus $T$ over $R_{\pi}$, such that $(A^0_{f/R_{\pi}})\hat{\ }\simeq \widehat{C}$. So we have an induced homomorphsim ${\mathcal A} \ra (A^0_{f/R_{\pi}})\h \ra \widehat{B}$, hence an induced morphism
 ${\mathcal X}_!\ra \widehat{B}$, hence an induced morphism
 $J^n_{\pi}({\mathcal X}_!)\ra J^n_{\pi}(B).$ Furthermore note that by the  N\'{e}ron property over $R_{\pi_0}$, plus base change from $R_{\pi_0}$ to $R_{\pi}$,  one has induced endomorphisms $T(l)_*$ of $A_{R_{\pi}}, A^0_{R_{\pi}}$, and ${\mathcal A}$ over $R_{\pi}$ for all primes $l$.
 Similarly, by the N\'{e}ron property over $R_{\pi}$, $\cO_f$ acts on $A_{f/R_{\pi}}$ and $A^0_{f/R_{\pi}}$. Clearly this action is compatible with the Hecke  action on $A^0_{R_{\pi}}$. Similarly, for $d\in (\bZ/p\bZ)^{\times}$, the diamond operators act on $\cX_!$ and $\cA$.
  By the functoriality of the Raynaud extension (cf. \cite{FC}, pp. 33-34) we have an induced action of $\cO_f$ on $C$ and hence on $B$; so we have a ring homomorphism $\iota:\cO_f\ra End(B/R_{\pi})$.

  Assume first $B\neq 0$.
  Then, for any prime $l$,  and any $d\in (\bZ/p\bZ)^{\times}$, we have   commutative diagrams
  \begin{equation}
  \label{cougar}
  \begin{array}{rcl}
  {\mathcal A} & \stackrel{T(l)_*}{\longrightarrow} & {\mathcal A}\\
  \downarrow & \  & \downarrow\\
  \widehat{B} & \stackrel{\iota(a_l)}{\longrightarrow} & \widehat{B}\end{array}
 \ \ \ \  \text{and}\ \ \ \ 
  \begin{array}{rcl}
  \cX_! & \stackrel{\langle d \rangle_p}{\longrightarrow} & \cX_!\\
  \downarrow & \  & \downarrow\\
  {\mathcal A} & \stackrel{\langle d \rangle_p}{\longrightarrow} & {\mathcal A}\\
  \downarrow & \  & \downarrow\\
  \widehat{B} & \stackrel{\iota(\epsilon_p(d))}{\longrightarrow} & \widehat{B}
  \end{array}\end{equation}
   On the other hand assume $B=0$; then $C=T$.  So $\cO_f$ acts on $T$ via a homomorphism $\iota^{\perp}:\cO_f \ra End(T/R_{\pi})$.
   Enlarging $R_{\pi}$ again and replacing the schemes (\ref{doarme}) by their pull-backs to this new $R_{\pi}$ we may (and will!)  assume $T_{K_{\pi}}$ is split over $K_{\pi}$.  So the identity component of its locally (!) finite type N\'{e}ron model has the form $D={\mathbb G}_{m/R_{\pi}}^{\nu}$ and, by the universality property of the N\'{e}ron model, we get a homomorphism $T \ra D$ over $R_{\pi}$. By the N\'{e}ron property over $R_{\pi}$ we have a homomorphism 
   $\iota^{\perp}:\cO_f \ra End(D/R_{\pi})$
   and corresponding diagrams:
   
   \begin{equation}
   \label{cheeta}\begin{array}{rcl}
  {\mathcal A} & \stackrel{T(l)_*}{\longrightarrow} & {\mathcal A}\\
  \downarrow & \  & \downarrow\\
  \widehat{D} & \stackrel{\iota^{\perp}(a_l)}{\longrightarrow} & \widehat{D}\end{array}
  \ \ \ \ \ \ \text{and} \ \ \ \ \ \ \ 
  \begin{array}{rcl}
  \cX_! & \stackrel{\langle d \rangle_p}{\longrightarrow} & \cX_!\\
  \downarrow & \  & \downarrow\\
  {\mathcal A} & \stackrel{\langle d \rangle_p}{\longrightarrow} & {\mathcal A}\\
  \downarrow & \  & \downarrow\\
  \widehat{D} & \stackrel{\iota^{\perp}(\epsilon_p(d))}{\longrightarrow} & \widehat{D}
  \end{array}\end{equation}

 \subsection{Lift of Igusa curve to characteristic zero \cite{over}}\label{zzzzz}
 Recall that $L$ denotes the line bundle on $X_1(N)_{R_p}$
such that the sections of the powers of $L$ identify with the modular forms of various weights on $\Gamma_1(N)$; cf. \cite{Gross} p. 450. Let $E_{p-1}\in H^0(X_1(N)_{R_p},L^{p-1})$ be the normalized Eisenstein form of weight $p-1$ and let $(ss)$ be the supersingular locus on $X_1(N)_{R_p}$, defined as the  locus of $E_{p-1}$.
Here  $E_{p-1}$ is normalized by the condition that its Fourier expansion has constant term $1$.
Take an open  covering 
\begin{equation}
\label{coverica}
X_1(N)=\bigcup X_1(N)_i\end{equation}
  such that  $L$ is trivial on each $X_1(N)_i$ and  let $x_i$ be a basis of $L$ on $X_1(N)_i$. Then $E_{p-1}=\varphi_i x_i^{p-1}$
where $\varphi_i \in \cO(X_1(N)_i)$. Set $x_i=u_{ij}x_j$, $u_{ij}\in \cO^{\times}(X_1(N)_{ij})$, $X_1(N)_{ij}=X_1(N)_i \cap X_1(N)_j$. Set $z_i=x_i^{-1}$ (in \cite{over} this $z_i$ was denoted by $t_i$) and consider the $R_{\pi}$-scheme
$X_1(N)_{!!}$ obtained by gluing the schemes 
$$X_1(N)_{!!i}:=Spec\ \cO(X_1(N)_{i,R_{\pi}})[z_i]/(z_i^{p-1}-\varphi_i)$$ via $z_i=u_{ij}^{-1}z_j$ where 
$$X_1(N)_{i,R_{\pi}}:=X_1(N)_i\otimes_{R_p} R_{\pi}.$$
 In the discussion below denote by an upper bar the functor $\otimes k$.
Then $\overline{X_1(N)_{!!}}$ is clearly birationally isomorphic (hence isomorphic)
to the Igusa curve $I$ (cf. \cite{Gross}, pp. 460, 461).
Now let $X \subset X_1(N)$ over $R_p$ be an affine open set disjoint from the supersingular locus. 
Consider the open  covering 
\begin{equation}
\label{cover}
X=\bigcup X_i\end{equation}
where $X_i=X_1(N)_i\cap X$. Consider the $R_{\pi}$-scheme
$X_{!!}\subset X_1(N)_{!!}$ obtained by gluing the schemes $X_{!!i}:=Spec\ \cO(X_{i,R_{\pi}})[z_i]/(z_i^{p-1}-\varphi_i)$ via $z_i=u_{ij}^{-1}z_j$ (where $X_{i,R_{\pi}}:=X_i\otimes_{R_p} R_{\pi}$). Note that $z_i^{p-1}-\varphi_i$ are monic polynomials whose derivatives are invertible in $\cO(X_{i,R_{\pi}})[z_i]/(z_i^{p-1}-\varphi_i)$
so $X_{i!!}$ is an \'{e}tale cover of $X_{i,R_{\pi}}$ and is the pull back of the latter
in $X_1(N)_{i!!}$.
 Note that the scheme $\overline{X_{!!}}=X_{!!} \otimes k$ is isomorphic to $\overline{\mathcal X}_!={\mathcal X}_!\otimes k$. Recall  from \cite{over} that the isomorphism $\overline{X}_{!!} \simeq \overline{\mathcal X}_!$  lifts uniquely to an isomorphism $(X_{!!})\h\simeq {\mathcal X}_!$. This follows immediately by applying
 the  standard Lemma \ref{hensel} below to $S:=\cO({\mathcal X}_i)$, ${\mathcal X}_i=\widehat{X}_i$,
 $T=\cO({\mathcal X}_{!i})$, ${\mathcal X}_{!i}=e^{-1}({\mathcal X}_i)$.
 \qed

\begin{lemma}
 \label{hensel}
 Let $S \ra T$ be a homomorphism of flat $\pi$-adically complete $R_{\pi}$-algebras, let $f\in S[z]$ be a monic polynomial and assume we have a homomorphism $\overline{\sigma}:\overline{S}[z]/(\overline{f}) \ra \overline{T}$ such that $df/dz$ is invertible in $\overline{S}[z]/(\overline{f})$. Then $\overline{\sigma}$ lifts uniquely to a homomorphism
$\sigma:S[z]/(f)\ra T$. If $\overline{\sigma}$ is an isomorphism then so is $\sigma$.
\end{lemma}

{\it Proof}. One defines $\sigma$ by sending the class $\xi\in S[z]/(f)$ of  $z$ into the unique solution $\tau\in T$  (which exists by Hensel) of the equation $f=0$ with $\overline{\sigma}(\overline{\xi})= \overline{\tau}$, where the upper bar means class mod $\pi$. \qed
  
   \subsection{Serre lifts}
 Let $M$ be any positive integer. 
In what follows a {\it classical modular form} over a ring $B$, of weight $\kappa$, on $\Gamma_1(M)$
 will be understood in the sense of \cite{Gross} as a rule that attaches to any $B$-algebra $C$
and any triple consisting of an elliptic curve  $E/C$, an embedding $\mu_{M,C} \ra
E[M]$, and an invertible one form on $E$ an element of $C$ satisfying the usual compatibility rules and the usual holomorphy condition for the Fourier expansion (evaluation on the Tate curve).
We denote by
$M_{\kappa}(\Gamma_1(M),B)$
the $B$-module of all these forms and by
$S_{\kappa}(\Gamma_1(M),B)$
the module of corresponding cusp forms. 
 For the concepts of  Hecke operators,  diamond operators, and newforms used below we refer to \cite{Gross}. Let  $\Theta_{p}:({\mathbb Z}/p{\mathbb Z})^{\times}\ra \bZ_p^{\times}\subset R_{\pi}^{\times}$ be the Teichm\"{u}ller character.
We have the following ``reduction to weight $2$" result \cite{Gross}, p. 478:

\begin{proposition}
\label{dudee}
Let $(N,p)=1, N>4$ and let $\overline{f}\in S_{\kappa}(\Gamma_1(N),k)$ be a newform, where $3\leq \kappa \leq p$. Then there exists   a newform $f\in S_2(\Gamma_1(Np),{\mathbb C})$ with character $\epsilon$ 
and there exists  an embedding 
$\rho:\cO_f\ra R_{\pi}$
such that:

1) The  image  of $f$ via $\cO_f \stackrel{\rho}{\ra}R_{\pi}\ra k$ is $\overline{f}$; 

2)  $\rho\circ \epsilon_p=\Theta_{p}^{\kappa-2}$.
\end{proposition}

In the above $\epsilon_p$ is the ``$p$-part" in the decomposition $\epsilon=\epsilon_p\epsilon_N$. We shall refer to the pair $(f,\rho)$, or simply to $f$, as a {\it Serre lift} of $\overline{f}$. Recall from \cite{Gross} that Serre lifts are not unique.

 \section{Construction of differential modular forms}

 \subsection{$\d_{\pi}$-characters}\label{dch}

 Let $B$ be an abelian scheme over $R_{\pi}$. Recall that by a $\d_{\pi}$-character of $B$, of order $n$, we understand
 a homomorphism $J^n_{\pi}(B) \ra \widehat{{\mathbb G}_a}$.
 
 \begin{proposition}
 Consider an abelian scheme $B$ over $R_{\pi}$ of relative dimension $g$. Then the $R_{\pi}$-module of $\d_{\pi}$-characters of $B$ of order $2$ has rank $\geq g$.
 \end{proposition}
 
 In \cite{char} this was proved under the assumption that the ramification index of $R_{\pi}$ satisfies $e\leq p-2$; the method there does not seem to extend to arbitrary $e$. The argument below, for general $e$, follows a method in \cite{fermat}.
 
 \medskip
 
 {\it Proof}. Let ${\mathcal F}\in R_{\pi}[[T_1,T_2]]^g$ be the formal group law of $B$ and $l(T)\in K_{\pi}[[T]]^g$ its logarithm; here $T,T_1,T_2$ are $g$-tuples of variables. For any $n\geq 1$ let
 $$L_{\pi}^n=\frac{1}{\pi} l(\phi^n(T))_{|T=0}\in K_{\pi}[[T,\d_{\pi}T,...,\d_{\pi}^n T]]$$
 where $\phi:K_{\pi}[[T,\d_{\pi}T,...]]\ra K_{\pi}[[T,\d_{\pi}T,...]]$, 
 $$\phi(T)=T^p+\pi\d_{\pi}T,\ \phi(\d_{\pi}T)=(\d_{\pi}T)^p+\pi\d_{\pi}^2T,...$$
  We claim that there is an integer $\nu\geq 1$ such that
 $$\pi^{\nu} L^n_{\pi} \in R_{\pi}[[T]][\d_{\pi}T,...,\d^n_{\pi}]\h;$$
 and, by the way, if in addition $e\leq p-1$ then one can take $\nu=0$. Indeed by \cite{haze}, p. 64, we have
 $$l(T)=\sum_{|\alpha|\geq 1} \frac{A_{\alpha}}{|\alpha|} T^{\alpha}$$
 for some $A_{\alpha}\in R_{\pi}$. Hence setting $G^n_{\pi}=\frac{1}{\pi}\phi^n(T)_{|T=0}$ we have
 $$L_{\pi}^n=\frac{1}{\pi}\sum_{|\alpha|\geq 1}\phi^n(A_{\alpha})\frac{\pi^{|\alpha|-1}}{|\alpha|}(G^n_{\pi})^{\alpha}$$
 Now if $v_{\pi}$ is the $\pi$-adic valuation we have
 $$v_{\pi}\left( \frac{\pi^{|\alpha|-1}}{|\alpha|}\right)\geq |\alpha|-1-e\cdot \log_p
 |\alpha|$$
 So the above goes to $\infty$ for $|\alpha|\ra \infty$; moreover the above is $\geq 0$ if $e\leq p-1$. This proves the claim. 
 Let $N^2$ be the kernel of $J^2_{\pi}(B)\ra \widehat{B}$.
Exactly as in \cite{fermat}
 the $g$-tuples $\pi^{\nu}L^2_{\pi}$ and $\pi^{\nu}L^1_{\pi}$ define $2g$ homomorphisms $N^2\ra \widehat{{\mathbb G}_a}$ that are $R_{\pi}$-linearly independent. So $Hom(N^2,\widehat{{\mathbb G}_a})$ has rank $\geq 2g$. We are done by considering the standard exact sequence (cf. \cite{char}):
 $$ 0=H^0(B,\cO)
 \ra Hom(J^2_{\pi}(B),\widehat{{\mathbb G}_a})\ra Hom(N^2,\widehat{{\mathbb G}_a})\ra H^1(B,\cO)\simeq R_{\pi}^g.$$
 \qed
 
 \begin{remark}
 Recall from \cite{char} that ${\mathbb G}_m=Spec R_{\pi}[x,x^{-1}]$ has a remarkable $\d_{\pi}$-character, i.e. a homomorphism $\psi:J^1_{\pi}({\mathbb G}_m)\ra \widehat{{\mathbb G}_a}$ defined by
 $$\psi = \pi^{m}\sum_{n\geq 1} (-1)^{n-1}\frac{\pi^n}{n}\left(\frac{\d_{\pi}x}{x^p}\right)^n$$
 for an appropriate $m\in \bZ$. So if $D={\mathbb G}_a^{\nu}$ then there is a canonical system
 $\psi_1,...,\psi_{\nu}:J^1_{\pi}(D)\ra \widehat{{\mathbb G}_a}$ of $\d_{\pi}$-characters obtained by composing $\psi$ above with the canonical projections onto the factors.
 \end{remark}
 
 \subsection{Semi-invariant $\d_{\pi}$-characters}
 Assume the notation of Section \ref{model}
  and let $\rho:K_f(\zeta_N,\zeta_p)\ra K_{\pi}$ be an embedding. (One can  always choose  $\pi$ such that such a $\rho$ exists.) 
 So $\rho(\cO_f[1/N,\zeta_N,\zeta_p])\subset R_{\pi}$. 
 Then, exactly as in \cite{eigen}, Proposition 4.5 and using the same argument as on the top of page 992, loc. cit., we get:
 
 \begin{lemma}
 \label{create}
 Assume $B\neq 0$ 
 and let ${\mathcal T}=\iota(\cO_f)\subset End(B/R_{\pi})$. Then there exists a non-zero
 $\d_{\pi}$-character $\psi\in \cO(J^2_{\pi}(B))$ and a ring homomorphism $\chi:{\mathcal T}\ra R_{\pi}$ such that 
 $\psi\circ \tau=\chi(\tau)\cdot \psi$
 for all $\tau\in {\mathcal T}$.
 \end{lemma}
 
 Similarly we get:

 \begin{lemma}
 \label{created}
 Assume $B=0$ (hence $D$ is defined and non-zero).  
 Let ${\mathcal T}=\iota^{\perp}(\cO_f)\subset End(D/R_{\pi})$. Then there exists a
 non-zero $\d_{\pi}$-character $\psi\in \cO(J^1_{\pi}(D))$ and a ring homomorphism $\chi:{\mathcal T}\ra R_{\pi}$ such that 
 $\psi\circ \tau=\chi(\tau)\cdot \psi$
 for all $\tau\in {\mathcal T}$.
 \end{lemma}

\subsection{$\d_{\pi}$-modular forms from $\cX_!$}

Recall from \cite{over} the following:

\begin{proposition}
There is a natural isomorphism
$$\cO(J^r_{\pi}(\cX_!)) \simeq \bigoplus_{\kappa=0}^{p-2}M_{\pi}^r(-\kappa).$$
Via this isomorphism $M_{\pi}^r(-\kappa)$ is identified with the subspace of 
all $\varphi \in \cO(J^r_{\pi}(\cX_!))$ on which the diamond operators $\langle d \rangle_p$ with  $d\in (\bZ/p\bZ)^{\times}$ act by the character $d\mapsto \Theta_p(d)^{-\kappa}$ (i.e. $\langle d \rangle_p \varphi=\Theta_p(d)^{-\kappa} \varphi$.)
\end{proposition}

The ring $\cO(J^n_{\pi}(\cX_!))$ can be referred to as the ring of {\it $\d_{\pi}$-modular forms} of level $Np$, order $n$, and weight $0$.

\bigskip

Let now $f$  be as in section \ref{model} and let $\rho:\cO_f[\zeta_{Np}]\ra R_{\pi}$ be an embedding. Assume the notation in section \ref{model} and let $\epsilon=\epsilon_N\epsilon_p$ be the character of $f$. 

Assume first $B\neq 0$.

 Let $\psi:J_{\pi}^2(B)\ra \widehat{{\mathbb G}_a}=\widehat{{\mathbb A}^1}$ and $\chi$ be as in Lemma 
 \ref{create}. 
 So for $d\in (\bZ/(p\bZ)^{\times}$ we have
  $$\psi \circ \iota(\epsilon_p(l))=
\chi (\iota(\epsilon_p(d)))\cdot \psi.$$ Then if $\beta:\cX_!\ra \widehat{B}$ is the vertical morphism
in equation \ref{cougar} and if $\beta^2:J^2_{\pi}(\cX_!)\ra J^2_{\pi}(\widehat{B})=J^2_{\pi}(B)$ is the induced morphism set
\begin{equation}
\label{crocodile}
f^{\sharp}=\psi \circ \beta^2:J^2_{\pi}(\cX_!)\ra \widehat{\mathbb A^1}\end{equation}
We get
$$\langle d \rangle_p f^{\sharp}=\psi \circ \beta^2 \circ \langle d \rangle_p=\psi \circ \iota(\epsilon_p(d))\circ \beta^2=\chi (\iota(\epsilon_p(d)))\cdot \psi \circ \beta^2=
\chi (\iota(\epsilon_p(d)))\cdot f^{\sharp}.$$
Set $\rho':=\chi\circ \iota:\cO_f\ra R_{\pi}$.
Then $f^{\sharp}$ has weight $-\kappa'$ for some integer $0\leq \kappa'\leq p-2$ if and only if
$$\rho' \circ \epsilon_p=\Theta_p^{-\kappa'}.$$
If in addition  $f$ and $\rho$ are as in Proposition \ref{dudee}, so in particular $\overline{f}$ has weight $\kappa$, then
$$\rho \circ \epsilon_p=\Theta_p^{\kappa-2}.$$
Let $g\in (\bZ/p\bZ)^{\times}$ be a generator; so $\Theta_p(g)\in \bZ_p^{\times}$ is a primitive $(p-1)$-th root of unity in ${\mathbb Q}_p$ and $\epsilon_p(g)$ is a (not necessarily primitive) $(p-1)$-th root of unity in ${\mathbb C}$.  Then $\rho,\rho'$ induce embeddings $r,r':{\mathbb Q}(\epsilon_p(g))\ra {\mathbb Q}(\Theta_p(g))$ hence  
$r(\epsilon_p(g))=\Theta_p(g)^{\kappa-2}$. Now since 
${\mathbb Q}(\Theta_p(g))/{\mathbb Q}$ is Galois
there is an automorphism $\gamma$ 
of  this extension such that $\gamma \circ r=r'$; by the structure of the Galois group
we must have $\gamma(\Theta_p(g))=\Theta_p(g)^c$ for some integer $c$ coprime to $p-1$. We get
$$\rho'(\epsilon_p(g))=r'(\epsilon_p(g)))=\gamma(r(\epsilon_p(g))=\Theta_p(g)^{c(\kappa-2)}$$
hence
$\rho'\circ \epsilon_p=\Theta_p^{c(\kappa-2)}=\Theta_p^{-\kappa'}$, where $\kappa'$ is the unique integer between $0$ and $p-2$ such that 
$\kappa'\equiv c(2-\kappa)$ mod $p-1$.
Then we have proved the following:

\begin{proposition}\label{11}
Let $\overline{f}\in S_{\kappa}(\Gamma_1(N),k)$ be a newform as in Proposition \ref{dudee} and let $(f,\rho)$ be a Serre lift (cf. loc.cit.). 
 Assume $B\neq 0$ and let $f^{\sharp}$ be as in \ref{crocodile}. Then  
 $f^{\sharp}\in M^2_{\pi}(-\kappa')$ where $\kappa'$ is a conjugate of $\kappa$.
 \end{proposition}

Similarly if $B=0$ (hence $D$ is defined and non-zero) we let $\psi:J_{\pi}^1(D)\ra \widehat{{\mathbb G}_a}=\widehat{{\mathbb A}^1}$ and $\chi$ be as in Lemma 
 \ref{created}. Then if $\beta:\cX_!\ra \widehat{D}$ is the vertical morphism
in equation \ref{cheeta} and if $\beta^1:J^1_{\pi}(\cX_!)\ra J^2(\widehat{D})=J^2_{\pi}(D)$ is the induced morphism we set
\begin{equation}
\label{frog}
f^{\sharp}=\psi \circ \beta^1:J^1_{\pi}(\cX_!)\ra \widehat{\mathbb A^1}\end{equation}
and we get:

\begin{proposition}\label{22}
Let $\overline{f}\in S_{\kappa}(\Gamma_1(N),k)$ be a newform as in Proposition \ref{dudee} and let $(f,\rho)$ be a Serre lift (cf. loc.cit.). 
  Assume $B=0$ and let $f^{\sharp}$ be as in \ref{frog}. Then  
 $f^{\sharp}\in M^1_{\pi}(-\kappa')$ where $\kappa'$ is a conjugate of $\kappa$.
 \end{proposition}
 
 The proof of Theorem \ref{main}
 will be concluded if we prove the following Proposition, where we assume that the embedding $\rho:\cO_f \ra R_{\pi}$ is extended (which is always possible by changing $\pi$) to an embedding of the Galois closure of $K_f$ into $K_{\pi}$.
 
 \begin{proposition}
 \label{33}
 Let $f^{\sharp}$ be as in \ref{crocodile} or \ref{frog}. Then $f^{\sharp}$ is non-zero and its $\d_{\pi}$-Fourier expansion has the form 
 $$E(f^{\sharp})=\sum_{\sigma}P_{\sigma}(\phi)(f^{(-1)}(q))^{\sigma},$$
with $P_{\sigma}(\phi)$   polynomials in $K_{\pi}[\phi]$ and $\sigma$ running through the set of all embeddings of $K_f$ into ${\mathbb C}$.
\end{proposition}
 
 {\it Proof}.  We have $f^{\sharp}\neq 0$ because $\psi\neq 0$ and the images of $\beta:\cX_!\ra \widehat{B}$ and  $\beta:\cX_!\ra \widehat{D}$ generate $\widehat{B}$ and $\widehat{D}$ respectively. The shape of $E(f^{\sharp})$  follows directly from the way $f^{\sharp}$ was defined, exactly as in \cite{eigen}, proof of Theorem 6.3.
For convenience we sketch here the argument. Assume  $B\neq 0$; the case $B=0$ is similar (with $D$ playing the role of $B$). Let $g$ be the dimension of $B$, let $\omega_1,...,\omega_g$ be a basis of global $1$-forms on $B$ and consider the formal logarithm $L=L(T_1,...,T_g)$ of $B$ with components $L_1,...,L_g$ (which are series in the variables $T_1,...,T_g$) such that $\omega_i=dL_i$. One can write 
$$\beta^*\omega_i=\sum_{\sigma}c_{i\sigma}\sum_n a_n^{\sigma}q^{n-1}dq,$$
with $c_{i\sigma}\in R_{\pi}$. On the other hand, if $\beta^*(T_i)=\varphi_i(q)\in R_{\pi}[[q]]$ and $\varphi=(\varphi_1,...,\varphi_g)$, we have
$$\beta^*\omega_i=\beta^*\left(\sum_j \frac{dL_i}{dT_j}dT_j\right)=\sum_j \frac{dL_i}{dT_j}(\varphi(q))\frac{d\varphi_j}{dq}dq=\frac{d}{dq} (L_i(\varphi(q))dq.$$
One then gets
$$L_i(\varphi(q))=
\sum_{\sigma}c_{i\sigma}\sum_n a_n^{\sigma} \frac{q^n}{n}
=\sum_{\sigma}c_{i\sigma}(f^{(-1)}(q))^{\sigma}.$$
On the other hand one has
$$\psi=\sum_i Q_i(\phi)L_i$$
for some polynomials $Q_i(\phi)\in K_{\pi}[\phi]$. We get
$$f^{\sharp}=\beta^*\psi=\sum_i Q_i(\phi)\left(\sum_{\sigma}c_{i\sigma}(f^{(-1)}(q))^{\sigma}\right)=\sum_{\sigma}\left(\sum_i c_{i\sigma}Q_i(\phi)\right) (f^{(-1)}(q))^{\sigma},$$
and we are done by setting $P_{\sigma}(\phi)=\sum_i c_{i\sigma}Q_i(\phi)$.
 \qed

\section{Link with Igusa differential modular functions}

The aim of this section is to construct a natural (and somewhat unexpected)  ring homomorphism from the ring $\cO(J^n_{\pi}(\cX_!))$ of $\d_{\pi}$-modular forms of level $Np$, order $n$, and weight $0$ to the ring $S^{n+1}_{\heartsuit}\otimes_{R_p}R_{\pi}$ where $S^{n+1}_{\heartsuit}$ is the ring of {\it Igusa  $\d_p$-modular forms of level $N$} \cite{igusa}. We start by reviewing the rings
$S^n_{\heartsuit}$.

Recall our basic setting in section 2.
With  $X\subset X_1(N)$ disjoint from the supersingular locus, and with $L$ and $V\ra X$ as in section 2 set:
$$S^n=\cO(J^n_p(X)),\ M^n=\cO(J^n_p(V)),\ S^{\infty}=\lim_{\ra} S^n,\ M^{\infty}=\lim_{\ra} M^n;$$ the elements of $M^n$ are the $\d_p$-modular functions of order $n$. 
For $w\in W$ recall that we let $M^n(w)$ be the the $R_p$-module of all $f\in M^n$ such that $\lambda \star f=\lambda^w f$ for $\lambda\in R_p^{\times}$, where $\star$ is the natural action of the multiplicative group; the elements of $M^n(w)$ are the $\d_p$-modular forms of weight $w$. On the other hand recall the natural $\d_p$-expansion map
$$E:M^n\ra R_p((q))^n:=R_p((q))[q',...,q^{(n)}]\h;$$
its restriction to each $M^n(w)$ is injective.
Let $S^n_{\heartsuit} \subset R_p((q))^n$ be the image of $M^n$ via $E$ and let $$S^{\infty}_{\heartsuit}=\lim_{\ra} S^n_{\heartsuit}\subset R_p((q))^{\infty}:=\lim_{\ra} R_p((q))^n.$$
The ring $S^{\infty}_{\heartsuit}$ can be  referred to as the ring of {\it Igusa $\d_p$-modular functions of level $N$} (or the ring of the $\d_p$-Igusa curve of level $N$ \cite{igusa}).
Recall from \cite{Barcau} and \cite{book}, p. 269,  that there is a (necessarily unique) $f^{\partial}\in M^1(\phi-1)$ such that
$E(f^{\partial})=1$. 
Moreover we have a congruence
\begin{equation}
\label{hasse}
f^{\partial}\equiv E_{p-1}\ \ \text{mod}\ \ p\ \ \text{in}\ \ M^1.
\end{equation}

\begin{theorem}\label{natur}
There exists a sequence of $S^0$-algebra homomorphisms,
$$\cO(J^n_{\pi}(\cX_!))\ra S^{n+1}_{\heartsuit}\otimes_{R_p} R_{\pi},$$
compatible with each other as $n$ varies, and commuting with $\d_{\pi}$.
\end{theorem}

As we shall see these homomorphisms are ``natural". Also if we consider the limit
$$\cO(J^{\infty}_{\pi}(\cX_!)):=\lim_{\stackrel{\ra}{n}} \cO(J^n_{\pi}(\cX_!))$$
 our Theorem \ref{natur} yields an $S^0$-algebra homomorphism, commuting with $\d_{\pi}$,
\begin{equation}
\label{natural2}
\cO(J^{\infty}_{\pi}(\cX_!))\ra S^{\infty}_{\heartsuit}\otimes_{R_p} R_{\pi}.\end{equation}

\medskip

{\it Proof of Theorem \ref{natur}}. We begin by constructing an $S^0$-algebra homomorphism
\begin{equation}
\label{natural}
\cO(\cX_!)\ra S^1_{\heartsuit}\otimes_{R_p} R_{\pi}.
\end{equation}
Indeed
consider a cover as in Equation \ref{cover} and, using the notation preceding that equation,  set $f^{\partial}=\Phi_i x_i^{\phi-1}$ where $\Phi_i\in S^1_i:=\cO(J^n_p(X_i))$. Let $$M^1_i=S^1_i[x_i,x_i^{-1},\d_p x_i]\h=S^1_i[z_i,z_i^{-1},\d_p z_i]\h$$ Then
$$M^1_i/(f^{\partial}-1)\simeq S^1_i[z_i,z_i^{-1},\d_p z_i]\h/(z_i^{\phi-1}-\Phi_i).$$
Since the latter is Noetherian and $p$-adically complete, and since
$$(M^1_i/(f^{\partial}-1))\otimes_{R_p}R_{\pi}$$ is finite over it, it follows that $(M^1_i/(f^{\partial}-1))\otimes_{R_p}R_{\pi}$ is $p$-adically complete, hence $\pi$-adically complete.
Recall that, setting $S_i=\cO(X_i)$, we have
$$\cO(\cX_{!i})=\frac{(S_i\otimes_{R_p}R_{\pi})[z_i]\h}{(z_i^{p-1}-\varphi_i)},\ \ (M^1_i/(f^{\partial}-1))\otimes_{R_p}R_{\pi}=\frac{S^1_i[z_i,z_i^{-1},\d_p z_i]\h \otimes_{R_p}R_{\pi}}{(z_i^{\phi-1}-\Phi_i)}.$$
Hence
$$\overline{\cO(\cX_{!i})}=\frac{\overline{S_i}[z_i]}{(z_i^{p-1}-\overline{\varphi_i})},\ \ \overline{(M^1_i/(f^{\partial}-1))\otimes_{R_p}R_{\pi}}=\frac{\overline{S^1_i}[z_i,z_i^{-1},\d_pz_i]}{(z_i^{p-1}-\overline{\Phi_i})}.$$
Since, by Equation \ref{hasse},  $\overline{\Phi_i}=\overline{\varphi_i}$ in $\overline{S^1_i}$ and $\overline{S_i}\subset \overline{S^1_i}$ we get  natural inclusions, compatible with varying $i$:
$$\overline{\cO(\cX_{!i})}\subset \overline{(M^1_i/(f^{\partial}-1))\otimes_{R_p}R_{\pi}}$$
and hence, by Lemma \ref{hensel}, unique liftings of the above to homomorphisms
$$\cO(\cX_{!i})\ra (M^1_i/(f^{\partial}-1))\otimes_{R_p}R_{\pi}.$$
By uniqueness these liftings glue together to give a homomorphism
$$\cO(\cX_{!})\ra (M^1/(f^{\partial}-1))\otimes_{R_p}R_{\pi}.$$
Composing with the natural map 
$$(M^1/(f^{\partial}-1))\otimes_{R_p}R_{\pi}\ra S^1_{\heartsuit}\otimes_{R_p}R_{\pi}$$
we get the desired homomorphism \ref{natural}

Now by
the universality property of $(\cO(J^n_{\pi}(\cX_!)))_{n\geq 0}$ and by the fact that 
$$(S^{n+1}_{\heartsuit}\otimes_{R_p} R_{\pi})_{n\geq 0}$$ is a $\d_{\pi}$-prolongation sequence we get that 
 the homomorphism \ref{natural} induces 
 natural compatible ring homomorphisms commuting with $\d_{\pi}$ as in the statement of the theorem.\qed

\begin{remark}
It would be interesting to see if the homomorphism \ref{natural2} is injective. On the other hand we claim that its reduction mod $\pi$,
\begin{equation}
\label{natural3}
\overline{\cO(J^{\infty}_{\pi}(\cX_{!}))}\ra \overline{S^{\infty}_{\heartsuit}},\end{equation}
is not injective. Indeed,
since $X_{!!}$ is \'{e}tale over $X$ it follows from \cite{char}, Proposition 1.4,  that 
$$\cO(J^n_{\pi}(\cX_{!i}))=\cO(\cX_{!i})[\d_{\pi}z_i,...,\d_{\pi}^nz_i]\h$$
hence
$$\overline{\cO(J^n_{\pi}(\cX_{!i}))}=\overline{\cO(\cX_{!i})}[\d_{\pi}z_i,...,\d_{\pi}^nz_i].$$
On the other hand by \cite{igusa}, Introduction,  we have 
$$\overline{S^{\infty}_{\heartsuit}\otimes_{R_p}R_{\pi}}=\overline{S^{\infty}_{\heartsuit}}=\frac{\overline{S^{\infty}}[z_i,\d_p z_i, \d_p^2 z_i,...]}{(f^{\partial}-1,\d_p(f^{\partial}-1),\d_p^2(f^{\partial}-1),...)}.$$
So the map
\ref{natural3}
 sends $\d^a_{\pi}z_i$ into $0$ for all $1\leq a \leq p-2$; in particular the map \ref{natural3} is not injective.
\end{remark}

\bibliographystyle{amsplain}

\end{document}